\documentclass{ifacconf}

\usepackage{natbib}

\usepackage{amsmath}
\usepackage{amssymb}
\usepackage{mathtools}
\usepackage{array}

\usepackage{booktabs}
\usepackage{tikz}
\usetikzlibrary{arrows, arrows.meta}
\usepackage{algorithmic}
\usepackage{graphicx}
\usepackage{textcomp}
\usepackage{xcolor}
\usepackage{enumitem}
\usepackage{url}

\DeclareMathOperator{\Real}{Re}
\DeclareMathOperator{\Imag}{Im}

\DeclarePairedDelimiter{\abs}{\lvert}{\rvert}

\newcommand{\suchthat}{\ifnum\currentgrouptype=16 \mathrel{}\middle|\mathrel{}\else\mid\fi}

\newlist{steps}{enumerate}{1}
\setlist[steps]{label={\arabic*.}, ref={\arabic*}, leftmargin=*}

\theoremstyle{definition}
\newtheorem{remark}{Remark}

\begin{document}

\begin{frontmatter}

\title{New Features of P3$\delta$ software: Partial Pole Placement via Delay Action}

\author[First,Second]{Islam Boussaada}
\author[First]{Guilherme Mazanti}
\author[First]{Silviu-Iulian Niculescu}
\author[Second]{Adrien Leclerc}
\author[Second]{Jayvir Raj}
\author[Second]{Max Perraudin}
\address[First]{Universit\'e Paris-Saclay, CNRS, CentraleSup\'elec, Inria, Laboratoire des Signaux et Syst\`emes, 91190, Gif-sur-Yvette, France. (e-mail: firstname.lastname@l2s.centralesupelec.fr).}
\address[Second]{IPSA, 63 Boulevard de Brandebourg, 94200 Ivry-sur-Seine, France}


\begin{abstract}
This paper presents the software \emph{Partial Pole Placement via Delay Action}, or P3$\delta$ for short.
P3$\delta$ is a Python software with a friendly user interface for the design of parametric stabilizing feedback laws with time-delays, thanks to two properties of the distribution of quasipolynomials' zeros, called multiplicity-induced-dominancy and coexisting real roots-induced-dominancy. After recalling recent theoretical results on these properties and their use for the feedback stabilization of control systems operating under time delays, the paper presents the main features of the current version of P3$\delta$. We detail, in particular, the assignable admissible region (the set of allowable dominant roots and the corresponding delay), which helps the user in the choice of input information, allowing a reliable stabilizing delayed feedback. We also present the newly set online version of P3$\delta$.
\end{abstract}
\end{frontmatter}

\section{Introduction}

Even though time delays commonly lead to desynchronizing or destabilizing effects on the dynamics of the system they appear, some works have emphasized that delays may also have stabilizing effects in control design. One of the earliest such works is \cite{Tallman1958Analog}, where authors use a delay-based controller to improve the stability of systems with oscillatory behavior and small damping. Another strategy, used for instance in \cite{Suh1979Proportional} and \cite{Atay1999Balancing}, is to replace the classical proportional-derivative controller by a proportional-delayed controller, using a so-called ``average derivative action'' due to the delay. It should also be noted that closed-loop stability may be guaranteed for some control systems precisely by the existence of the delay, a fact highlighted in particular in \cite{Niculescu2010Delay}, in which the reader may find further discussions on the stabilizing effects of delays. A growing literature exhibits the design of delayed controllers in a wide range of applications, such as the control of flexible mechanical structures or the regulation of networks (see, e.g., \cite{Boussaada2018Further, Irofti2016Codimension}).

In this paper, we present a software, P3$\delta$ (which stands for \emph{Partial Pole Placement via Delay Action}), which helps its user in the stability analysis and the stabilization of linear time-invariant delay-differential equations of retarded or neutral type with a single time delay, under the form
\begin{multline}
\label{MainEqn}
y^{(n)}(t) + a_{n-1} y^{(n-1)}(t) + \dotsb + a_0 y(t) \\ + b_m y^{(m)}(t - \tau) + \dotsb + b_0 y(t - \tau) = 0,
\end{multline}
where $\tau > 0$ is the positive delay, $y$ is the real-valued unknown function, $n$ and $m$ are nonnegative integers with $n \geq m$, and $a_0, \dotsc, a_{n-1}, b_0, \dotsc, b_m$ are real coefficients. Equation \eqref{MainEqn} is said to be of \emph{retarded type} if $n>m$, that is, the highest-order derivative only appears in the non-delayed term $y^{(n)}(t)$, and it is said to be of \emph{neutral type} if $m=n$ and $b_n \neq 0$, which corresponds to a highest-order derivative appearing in both non-delayed term $y^{(n)}(t)$ and delayed one $y^{(n)}(t-\tau)$.

A typical situation in which equations under the form \eqref{MainEqn} arise is the delayed feedback stabilization of a linear time-invariant controlled differential equation of the form $y^{(n)} + a_{n-1} y^{(n-1)}(t) + \dotsb + a_0 y(t) = u(t)$, when applying a delayed feedback control $u(t) = - b_m y^{(m)}(t - \tau) - \dotsb - b_0 y(t - \tau)$. In that case, the choice of the free parameters $b_0, \dotsc, b_m$ in the feedback control will have an important influence on the behavior of the closed-loop system. We refer the reader to \cite{Boussaada2018Dominancy, Boussaada2020Multiplicity, Boussaada2018Further, MazantiMultiplicity, MazantiEffects} for more examples of systems under the form \eqref{MainEqn} to which the methods of the present paper can be applied.

A classical technique to address the stability analysis of linear time-invariant time-delay systems is by the use of spectral methods (see, e.g., \cite{Hale1993Introduction, Michiels2014Stability}), which consist on the study of the complex roots of a \emph{characteristic function} of the system. The characteristic function of \eqref{MainEqn} is
\begin{equation}
\label{Delta}
\textstyle \Delta(s) = s^n + \sum_{k=0}^{n-1} a_k s^k + e^{-s \tau} \sum_{k=0}^m b_k s^k,
\end{equation}
and \eqref{MainEqn} is exponentially stable if and only if the \emph{spectral abscissa} $\gamma = \sup\{\Real s \suchthat \Delta(s) = 0\}$ satisfies $\gamma < 0$. 

Characteristic functions of time-delay systems are \emph{quasipolynomials}, i.e., functions which can be written as a finite sum of polynomials multiplied by exponentials. Due to their applications in the spectral analysis of time-delay systems, the study of quasipolynomials has been the subject of several works, such as \cite{Berenstein1995Complex, Hale1993Introduction, Stepan1989Retarded, Wielonsky2001Rolle}. Except for the particular case where no exponentials are present and the quasipolynomial reduces to a polynomial, quasipolynomials have infinitely many roots in the complex plane. In the context of stabilization of a control system by a delayed feedback law, one only disposes of finitely many parameters in the feedback law to choose the location of these infinitely many roots and place them in order to guarantee a negative spectral abscissa, and hence exponential stability of the closed-loop system.

\newcommand{\DominancyRefs}{\cite{Amrane2018Qualitative, BedouheneReal, Boussaada2018Dominancy, Boussaada2020Multiplicity, Boussaada2018Further, MazantiMultiplicity, Mazanti2020Qualitative, Mazanti2020Spectral}}

A possible strategy to stabilize a time-delay system is to select the free parameters of the system in order to choose the location of finitely many roots while also guaranteeing that the \emph{dominant root}, i.e., the rightmost root on the complex plane, is among the chosen ones. This has been the subject of several recent works, such as \DominancyRefs{}. Contrarily to the strategy of finite spectrum assignment used, e.g., in \cite{Manitius1979Finite}, the controllers designed using these techniques do not render the closed-loop system finite-dimensional, but control instead its rightmost spectral value. These methods also extend to some partial differential equations, as detailed, for instance, in \cite{MazantiEffects}. 

Notice that there exist other pole placement paradigms for time-delay systems, such as the \emph{continuous pole-placement} introduced in \cite{Michiels2002}. Based on the continuous dependence of the characteristic roots on the controller parameters, this technique consists in shifting the unstable characteristic roots from $\mathbb{C}_+$ to $\mathbb{C}_-$ in a ``quasi-continuous'' way, subject to the constraint that, during this shifting action, stable characteristic roots are not crossing the imaginary axis from $\mathbb{C}_-$ to $\mathbb{C}_+$ (see also \cite{RAMLAA}).

Two main strategies have been used in the works \DominancyRefs{} to assign finitely many roots while guaranteeing that the rightmost root is among them. The first one consists in assigning a real root of maximal multiplicity and proving that this root is necessarily the rightmost root of the characteristic quasipolynomial, a property which has been named \emph{multiplicity-induced-dominancy}, or MID for short. The second strategy consists in assigning a certain amount of real roots, typically equally spaced for simplicity, and proving that the rightmost root among them is also the rightmost root of the characteristic quasipolynomial, a property which has been named \emph{coexisting real roots-induced-dominancy}, or CRRID for short.

The MID property for \eqref{MainEqn} was shown, for instance, in \cite{Boussaada2018Further} in the case $n = 2$ and $m = 0$, in \cite{Boussaada2020Multiplicity} in the case $n = 2$ and $m = 1$ (see also \cite{Boussaada2018Dominancy}), and in \cite{MazantiMultiplicity} in the case of any positive integer $n$ and $m = n-1$ (see also \cite{Mazanti2020Qualitative}). It was also studied for neutral systems of orders $1$ and $2$ in \cite{Ma,Benarab2020MID, MazantiEffects}, and extended to complex conjugate roots of maximal multiplicity in \cite{Mazanti2020Spectral}. The CRRID property was shown, for instance, in \cite{Amrane2018Qualitative} in the cases $(n, m) = (2, 0)$ and $(n, m) = (1, 0)$, and in \cite{BedouheneReal} in the case of any positive integer $n$ and $m = 0$.

In all the above cases, the maximal multiplicity of a real root or, equivalently, the maximal number of coexisting simple real roots is the integer $n + m + 1$. Furthermore, the idea to exploit the nature of (real or complex) open-loop roots in control design was proposed for second-order systems in \cite{Boussaada2020Multiplicity} and further extended for arbitrary order systems with real-rooted plants in \cite{Balogh20,Balogh21}.

Based on the results from \DominancyRefs{}, the Python software P3$\delta$ allows for the parametric design of stabilizing feedback laws with time delays using the MID and CRRID properties. It has the advantage to allow constructive methods more appropriate for understanding the effect of uncertainties on the spectrum distribution. The first version of P3$\delta$ was formerly described in \cite{BoussaadaP3D} and it covered the design of feedback laws for linear time-invariant differential equations with a single time delay under the form \eqref{MainEqn} using MID techniques in retarded case. The newer version of P3$\delta$, described in the present paper, benefits also from the CRRID property and treats both retarded and neutral equations.

Stability, robustness, or bifurcation aspects of time-delay systems have also been the subject of other recently developed softwares, such as YALTA \citep{Avanessoff2014Hinfty}, which performs $H_\infty$ stability analysis of time-delay systems with commensurate delays, TRACE-DDE \citep{Breda2009Trace}, for the computation of characteristic roots and stability charts of linear autonomous time-delay systems, DDE-BIFTOOL \citep{Engelborghs2002Numerical}, interested in the computation, continuation, and stability analysis of steady-state solutions of time-delay systems and their bifurcations, and QPmR \citep{Vyhlidal2014QPmR}, for the computation of roots of quasipolynomials. One of the major novelties of P3$\delta$ lies in addressing the stabilization of control systems with time delays by using the MID and CRRID properties to design stabilizing feedback laws, making use of both symbolic and numeric computations.

\section{Description of P3$\delta$}
\label{SecP3Delta}

P3$\delta$ is freely available for download on \url{https://cutt.ly/p3delta}, where installation instructions, video demonstrations, and the user guide are also available. Interested readers may also contact directly any of the authors of the paper.

Three modes are implemented in the current version of P3$\delta$: ``Generic MID'', ``Generic CRRID'', and ``Control-oriented MID''. In the generic modes, all coefficients of \eqref{MainEqn} are assumed to be available for choice, whereas, in the control-oriented mode, P3$\delta$ assumes that $a_0, \dotsc, a_{n-1}$ are fixed and only $b_0, \dotsc, b_m$ are free. These three modes are detailed in Sections~\ref{SecGenericMID}, \ref{SecGenericCRRID}, and \ref{SecControlMID} below.

\subsection{Generic MID mode}
\label{SecGenericMID}

The ``Generic MID'' mode was already implemented in the first version of P3$\delta$ and is described in \cite{BoussaadaP3D}. For sake of completeness, we also provide a description of this mode here.

In the ``Generic MID'' mode, the user inputs the values of the delay $\tau$ and of the desired real root $s_0$ and P3$\delta$ computes all coefficients $a_0, \dotsc, a_{n-1}, b_0, \dotsc, b_{m}$ ensuring that the value $s_0$ is a dominant root of $\Delta$ of maximal multiplicity $n + m + 1$, using the procedure from \cite{MazantiMultiplicity, BoussaadaGeneric}. To use the ``Generic MID'' mode, the user should proceed as follows:
\begin{steps}
\item Enter the values of the integers $n$ and $m$ appearing in the differential equation \eqref{MainEqn}.
\item\label{GenericMIDStepSelection} Select the ``Generic MID'' option in the drop-down menu ``--- Choose type ---''.
\end{steps}
After this selection, the window of the program is filled with places for the other inputs and outputs of P3$\delta$.
\begin{steps}[resume]
\item\label{GenericMIDStepInputs} Enter the values of the desired real root of maximal multiplicity $s_0$ and of the delay $\tau$ in the corresponding fields that appear below the drop-down menu.
\item\label{GenericMIDStepConfirm} Enter the bounds $x_{\min}, x_{\max}, y_{\min}, y_{\max}$ of the rectangle $[x_{\min}, x_{\max}] \times [y_{\min}, y_{\max}] \subset \mathbb C$ in which P3$\delta$ will look for roots of \eqref{Delta} and press the ``Confirm'' button.
\end{steps}
Once the ``Confirm'' button is pressed, P3$\delta$ computes and displays the values of the coefficients $a_0, \dotsc,\allowbreak a_{n-1},\allowbreak b_0, \dotsc,\allowbreak b_m$ ensuring that $s_0$ is a root of maximal multiplicity of the quasipolynomial $\Delta$ from \eqref{Delta}, performs a numerical computation of all roots of $\Delta$ within the selected rectangle, and plot these roots in the plot ``Roots'' at the lower left corner of the window. The latter computation is carried out using Python's \texttt{cxroots} module (see \cite{cxroots}).

Optionally, after the previous computations are completed, the user may also simulate some trajectories of the system in time domain. This can be done, after completing step \ref{GenericMIDStepConfirm} above, by the following steps:
\begin{steps}[resume]
\item\label{GenericMIDStepTimeSimulationBegin} Choose the type of the initial condition from the drop-down menu ``--- Initial Solution ---''.
\end{steps}
The currently supported types are ``Constant'', ``Polynomial'', ``Exponential'', and ``Trigonometric'', which corresponds to initial conditions of the forms $x(t) = c$, $x(t) = \sum_{k=0}^r c_k t^k$, $x(t) = A e^{\gamma t}$, and $x(t) = A \sin(\omega t + \varphi)$, respectively, where $c, r, c_0, \dotsc, c_r, A, \gamma, \omega, \varphi$ are constants to be chosen by the user and the initial condition is defined in the time interval $[-\tau, 0]$.
\begin{steps}[resume]
\item\label{GenericMIDStepT} Enter the simulation time $T$ in the corresponding box.
\item Enter the values of the constants appearing in the expression of the initial condition in the corresponding input boxes.
\item\label{GenericMIDStepTimeSimulationEnd} After entering all the constants, press ``Enter'' on the keyboard or click on the ``Confirm'' button appearing in the same frame as the constants.
\end{steps}
After these steps, the numerical solution corresponding to the chosen initial condition is computed using an explicit Euler scheme in the time interval $[-\tau, T]$ and plotted in the graph on the ``Solutions'' part of the screen.

\subsection{Generic CRRID mode}
\label{SecGenericCRRID}

In the newly implemented ``Generic CRRID'' mode, the user inputs the values of the delay $\tau$ and of $n + m + 1$ desired real roots $s_1\geq\dotsb\geq s_{m+n+1}$ and P3$\delta$ computes all coefficients $a_0, \dotsc, a_{n-1}, b_0, \dotsc, b_{m}$ ensuring that the value $s_1$ is simple and dominant root of $\Delta$, using the procedure described in \cite{BedouheneReal}.

The use of the ``Generic CRRID'' mode is very similar to that of the ``Generic MID'' mode described in Section~\ref{SecGenericMID}, with the differences that, in Step~\ref{GenericMIDStepSelection}, the user should select the ``Generic CRRID'' mode and, in Step~\ref{GenericMIDStepInputs}, the user should enter the values of the desired real roots $s_1\geq\dotsb\geq s_{m+n+1}$ and of the delay $\tau$ in the corresponding fields.

As in the ``Generic MID'' mode, once the ``Confirm'' button is pressed, P3$\delta$ will compute and show the values of the coefficients $a_0, \dotsc, a_{n-1}, b_0, \dotsc, b_m$ ensuring that the real roots $s_1\geq\dotsb\geq s_{m+n+1}$ are simple roots of the quasipolynomial $\Delta$ from \eqref{Delta}, and numerically compute all roots of $\Delta$ within the selected rectangle using Python's \texttt{cxroots} module. The user can simulate trajectories in time domain by following Steps~\ref{GenericMIDStepTimeSimulationBegin}--\ref{GenericMIDStepTimeSimulationEnd} from Section~\ref{SecGenericMID}.

\subsection{Control-oriented MID mode}
\label{SecControlMID}

Contrarily to the previous two modes, the ``Control-oriented MID'' mode considers that the coefficients $a_0, \dotsc,\allowbreak a_{n-1}$ corresponding to the non-delayed terms of \eqref{MainEqn} are given and that the coefficients $b_0, \dotsc, b_m$ corresponding to the delayed terms are available for choice. As described in \cite{Boussaada2020Multiplicity, BoussaadaP3D}, in this mode, one imposes the existence of a real root of multiplicity $m + 2$ at $s_0$ by requiring that $\Delta(s_0) = \Delta^\prime(s_0) = \dotsb = \Delta^{(m+1)}(s_0) = 0$. This gives a system of $m + 2$ equations on the $m + 1$ unknowns $b_0, \dotsc, b_m$, and thus another parameter of the system, either the delay $\tau$ or the root $s_0$, must be considered as an unknown of the problem, while the value of the other is assumed to be fixed. Hence, in ``Control-oriented MID'' mode of P3$\delta$, the user must choose to input either the value of $\tau$ or that of $s_0$ (but not both) and P3$\delta$ computes all coefficients $b_0, \dotsc, b_m$ ensuring the existence of a dominant root of the quasipolynomial $\Delta$ from \eqref{Delta} of multiplicity $m+2$, as well as the value of the parameter $\tau$ or $s_0$ that has not been fixed by the user.

To use the ``Control-oriented MID'' mode, the user should proceed as in Section~\ref{SecGenericMID}, but selecting ``Control-oriented MID'' in Step~\ref{GenericMIDStepSelection} and then entering either $s_0$ or $\tau$, as well as the values of the known coefficients $a_0, \dotsc, a_{n-1}$. As in the previous two modes, the user may also select to perform a time-domain simulation by following the same steps.

\begin{remark}
In  the ``Control-oriented MID'' mode, it may happen to be impossible to choose a real root $s_0$ of multiplicity $m+2$. In this case, the previous version of P3$\delta$ warned the user of this fact and provided an equation relating $s_0$ and $\tau$. The user should either enter a value of $s_0$ such that this equation admits a positive root $\tau$ or a positive value of $\tau$ such that this equation admits a real root $s_0$ in order to proceed with the computations. Instead of providing this equation, the new version of P3$\delta$ provides the plot in the plane $(s_0, \tau)$ of the points for which the equation is satisfied. This region of admissible assignment should help the user in the choice of an allowable $s_0$ or $\tau$. This new feature is further described in Section~\ref{AdmissibilityRegi}.
\end{remark}

In the case where the user inputs the delay $\tau$ instead of the desired root $s_0$, the ``Control-oriented MID'' mode can also perform a numerical sensitivity analysis of the computed roots with respect to variations of $\tau$. The steps to get the sensitivity plot are the following:
\begin{steps}
\item Select the ``Sensitivity'' tab in the ``Roots'' plot.
\item Select ``tau sensitivity'' in the drop-down menu ``--- Sensitivity ---'' above the ``Roots'' plot.
\item Enter the value of the step $\varepsilon$ and the number of iterations $K$ in the corresponding boxes.
\item Enter the bounds $x_{\min}, x_{\max}, y_{\min}, y_{\max}$ of the rectangle $[x_{\min}, x_{\max}] \times [y_{\min}, y_{\max}] \subset \mathbb C$ in which P3$\delta$ will look for roots of \eqref{Delta}.
\end{steps}
Since the sensitivity computation may take quite some time, it is highly recommended to choose a smaller rectangle containing only a few roots of $\Delta$.
\begin{steps}[resume]
\item Press the ``Confirm'' button in the frame of the bounds of the rectangle.
\end{steps}
The sensitivity plot appears in the ``Roots'' plot and contains the roots of $\Delta$ in the selected rectangle for the values of delays $\tau + k \varepsilon$ for $k \in \{-K, -K+1, \dotsc, K-1, K\}$. Roots computed with negative values of $k$, corresponding to values of the delay smaller than $\tau$, are represented in shades of blue, with darker blue representing $k = -K$ and lighter tones representing increasing values of $k$. Roots computed with positive values of $k$, corresponding to values of the delay larger than $\tau$, are represented in shades of orange to red, with darker red representing $k = K$ and lighter tones moving to orange representing decreasing values of $k$. The roots computed with $k = 0$, corresponding to the nominal value of $\tau$ selected by the user, are represented by black diamonds.

\subsection{Assignment admissibility region}
\label{AdmissibilityRegi}

This section describes the admissibility region and how to display it using P3$\delta$. Note that this new feature of P3$\delta$ is available only in the ``Control-Oriented MID'' mode, since there are no constraints on $(s_0, \tau)$ in the other modes.

For given coefficients $a_0, \dotsc, a_{n-1}$, the admissibility region is defined as the set of pairs $(s_0, \tau) \in \mathbb R \times (0, +\infty)$ for which there exist real coefficients $b_0, \dotsc, b_m$ such that $s_0$ is a root of $\Delta$ of multiplicity at least $m + 2$ when the delay is $\tau$.

Let us first describe how this admissibility region can be determined. Consider the quasipolynomial $\Delta$ from \eqref{Delta} with known values for the coefficients $(a_i)_{0\leq i \leq n-1}$. We impose that $s_0 \in \mathbb R$ is a root of multiplicity at least $m + 2$ of $\Delta$. We use the $m+1$ equations $\Delta^{(k)}(s_0) = 0$, $k \in \{0, \dotsc, m\}$, in order to express the coefficients $(b_i)_{0 \leq i \leq m}$ in terms of $(a_i)_{0\leq i \leq n-1}$, $s_0$, and $\tau$. This is always possible, since these $m+1$ equations are linear in the $m+1$ variables $b_0, \dotsc, b_m$ and this linear system can be shown to admit a unique solution.

In order for $s_0$ to be a root of multiplicity $m + 2$ of $\Delta$, in addition to the $m+1$ equations $\Delta^{(k)}(s_0) = 0$, $k \in \{0, \dotsc, m\}$, the $m+1$-th derivative of $\Delta$ must also be zero at $s_0$, i.e., the equation $\Delta^{(m+1)}(s_0) = 0$ must also be satisfied. By replacing the previously found expressions of $(b_i)_{0 \leq i \leq m}$ into the equation $\Delta^{(m+1)}(s_0) = 0$, one obtains a relation between $s_0$, $\tau$, and the coefficients $(a_i)_{0 \leq i \leq n-1}$, which is a necessary and sufficient condition for $s_0$ to be a root of multiplicity at least $m+2$ of $\Delta$. Since the coefficients $(a_i)_{0 \leq i \leq n-1}$ are assumed to be known, the only unknowns in this equation are $s_0$ and $\tau$, and hence the admissibility region is the set of pairs $(s_0, \tau) \in \mathbb R \times (0, +\infty)$ satisfying this equation.

These computations are implemented symbolically in P3$\delta$ using Python's \texttt{sympy} module and, once the equation describing the admissibility region is obtained, the admissibility region is displayed to the user. Since it is not possible to display the full admissibility region in $\mathbb R \times (0, +\infty)$, the user is prompted for values $s_{0, \min} < 0$ and $\tau_{\max} > 0$ and only the part of the admissibility region inside the rectangle $[s_{0, \min}, 0] \times [0, \tau_{\max}]$ is displayed.

\begin{figure}[!h]
\centering
\includegraphics[width=0.8\columnwidth]{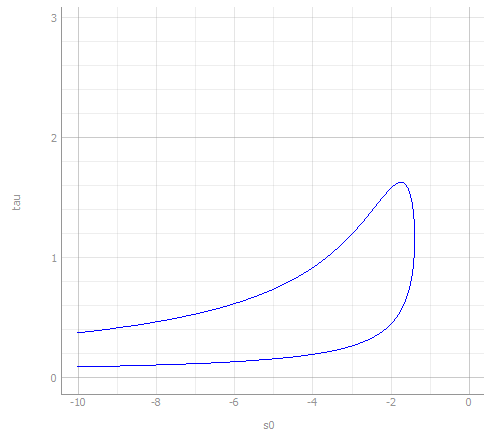}
\caption{Admissibility plot for a second-order retarded system with $a_0 = 1$, $a_1 = 1$, $\tau_{\max}=3$, and $s_{0, \min} = -10$.}
\label{FigAdmissibility}
\end{figure}

The admissibility region is always displayed when the user chooses the ``Control-oriented MID'' mode of P3$\delta$, after performing the following steps:
\begin{steps}
\item Enter the values of the integers $n$ and $m$ appearing in the differential equation \eqref{MainEqn}.
\item Select the ``Control-oriented MID'' option in the drop-down menu ``--- Choose type ---''.
\end{steps}
After this selection, the window of the program is filled with places for the other inputs and outputs of P3$\delta$.
\begin{steps}[resume]
\item Enter the values of the known coefficients $a_0, \dotsc, a_{n-1}$.
\item Enter the values of the limits $\tau_{\max} > 0$ and $s_{0, \min} < 0$ of the rectangle $[s_{0, \min}, 0] \times [0, \tau_{\max}]$ in which to plot the admissibility region in the respective fields ``tau limit'' and  ``s0 limit''.
\item Click on ``Confirm''.
\end{steps}
A new window appears with the admissibility region plotted. As an example, Figure~\ref{FigAdmissibility} shows the admissibility region in the case $n = 2$, $m = 1$, $a_0 = a_1 = 1$, and with $\tau_{\max} = 3$ and $s_{0, \min} = -10$.

After the admissibility region is shown, the user is prompted to enter a value for either $s_0$ or $\tau$ in order to proceed with the ``Control-oriented MID'' mode, as described in Section~\ref{SecControlMID}. The value of $s_0$ or $\tau$ entered by the user should be such that a point with the corresponding value of $s_0$ or $\tau$ exists in the admissibility region. In the example of Figure~\ref{FigAdmissibility}, the user cannot enter a value of $\tau$ exceeding $\approx 1.6$, nor a value of $s_0$ exceeding $\approx -1.5$.

\section{Illustrative example: Stabilization of an oscillator}
\label{SecExpl}

\begin{figure*}[ht]
\centering
\begin{tabular}{@{} >{\centering} m{0.33\textwidth} @{} >{\centering} m{0.34\textwidth} @{} >{\centering} m{0.33\textwidth} @{}}
\includegraphics[width=0.33\textwidth]{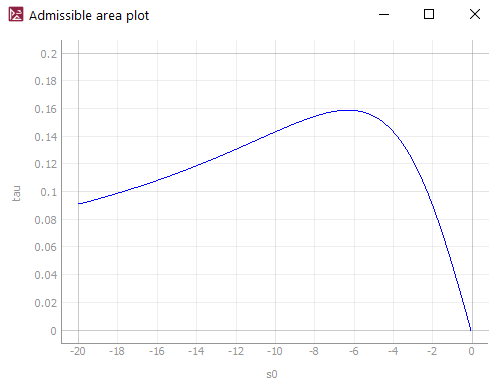} & \includegraphics[width=0.33\textwidth]{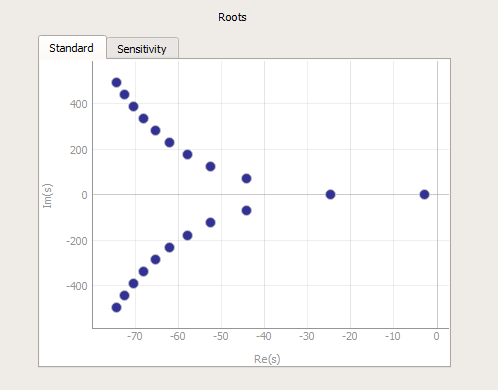} & \includegraphics[width=0.33\textwidth]{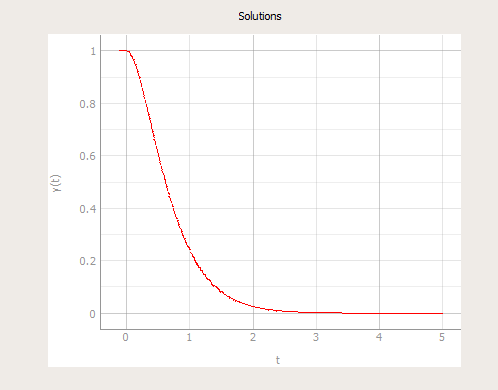} \tabularnewline
(a) & (b) & (c)
\end{tabular}
\caption{Example of the usage of P3$\delta$ in the ``Control-oriented MID'' mode for the delay-differential equation \eqref{EqExpl2} with $\omega = 2 \pi$. (a) Admissible region. (b) Roots of the characteristic quasipolynomial. (b) Numerical simulation of a solution with initial condition $x(t) = 1$ for $t \leq 0$.}
\label{FigExpl}
\end{figure*}

As an illustration of the use of P3$\delta$, let us consider the stabilization of an oscillator described by the equation $y^{\prime\prime}(t) + \omega^2 y(t) = u(t)$ with a delayed feedback control $u(t) = - b_0 y(t - \tau)$, which yields the delay-differential equation
\begin{equation}
\label{EqExpl2}
y^{\prime\prime}(t) + \omega^2 y(t) + b_0 y(t - \tau) = 0,
\end{equation}
whose characteristic quasipolynomial is $\Delta(s) = s^2 + \omega^2 + b_0 e^{-s \tau}$. This corresponds to \eqref{MainEqn} with $n = 2$, $m = 0$, $a_1 = 0$, and $a_0 = \omega^2$. Immediate computations show that $s_0 \in \mathbb R$ is a root of multiplicity at least $m+2 = 2$ if and only if
\begin{equation}
\label{Expl2Conditions}
b_0 = -e^{s_0 \tau} (s_0^2 + \omega^2), \qquad 2 s_0 + \tau(s_0^2 + \omega^2) = 0.
\end{equation}
The admissibility region is the set $\{(s_0, \tau) \in \mathbb R \times (0, +\infty) \suchthat 2 s_0 + \tau(s_0^2 + \omega^2) = 0\}$. For $\omega = 2 \pi$, we input $n = 2$ and $m = 0$ in P3$\delta$, select ``Control-oriented MID'', and choose $a_1 = 0$ and $a_0 = (2\pi)^2 \approx 39.48$. Selecting $s_{0, \min} = -20$ and $\tau_{\max} = 0.2$ for the admissibility plot, we get the admissibility region in Figure~\ref{FigExpl}(a). We then choose $\tau = 0.12$ in P3$\delta$ and obtain that $b_0 \approx -33.81$ and $s_0 \approx -2.859$, which agrees with \eqref{Expl2Conditions}. We also obtain the roots of $\Delta$ in a given rectangle, represented in Figure~\ref{FigExpl}(b) for the rectangle $\{s \in \mathbb C \suchthat \abs{\Real s} \leq 500,\, \abs{\Imag s} \leq 500\}$, and time simulations of solutions, for instance the one in Figure~\ref{FigExpl}(c) obtained with the constant initial condition $x(t) = 1$.

\section{Online version}

Since its creation, P3$\delta$ had vocation to be available to the greatest number and on all possible platforms. The current version of the software is available in local executable version, and now the development team wants to propose an online version ready to use in one click. The online version of P3$\delta$ is hosted on servers thanks to the \emph{Binder} service \citep{project_jupyter-proc-scipy-2018}. \emph{Binder} allows to create instances of personalized computing environment directly from a \emph{GitHub} repository that can be shared and used by users. The \emph{Binder} service is free to use and is powered by \emph{BinderHub}, an open-source tool that deploys the service in the cloud. The online version of P3$\delta$ is written in Python and structured as a \emph{Jupyter Notebook}, an open document format which can contain live code, equations, visualizations, and text. The Jupyter Notebook is completed by a friendly user interface built using interactive widgets from Python's \texttt{ipywidgets} module.

P3$\delta$ online is based on the program and features of the executable version. The first version of the online software includes for the moment features from the ``Generic MID'', ``Control-oriented MID'', and ``Generic CRRID'' modes of P3$\delta$ described above in Sections~\ref{SecGenericMID}--\ref{SecControlMID}, with user inputs similar to those from the executable version.

\begin{figure*}[ht]
\centering
\begin{minipage}[b]{0.333\textwidth}
\centering

\includegraphics[width=\textwidth]{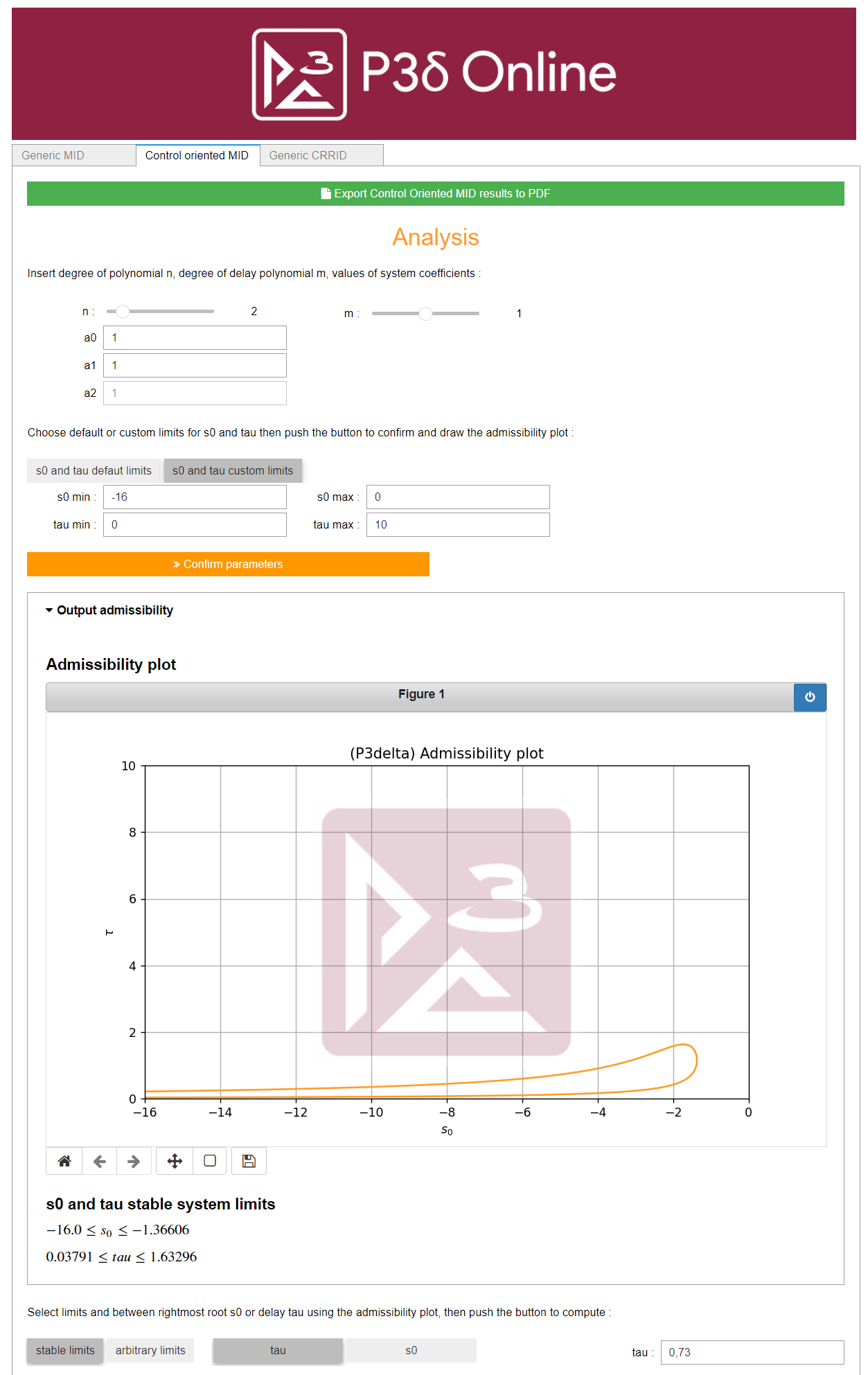}

(a)

\end{minipage}\begin{minipage}[b]{0.333\textwidth}
\centering

\includegraphics[width=\textwidth]{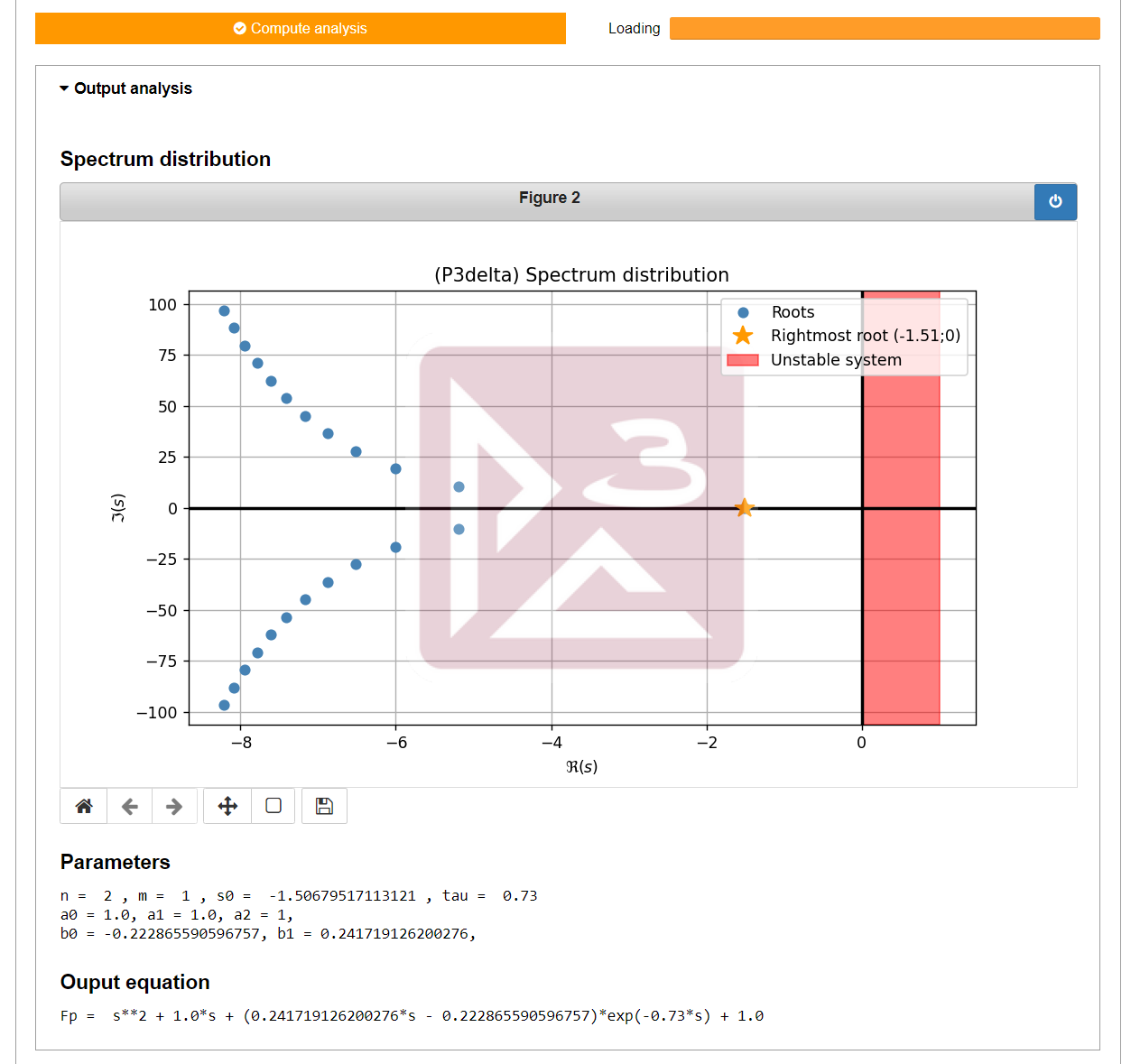}

(b)

\end{minipage}\begin{minipage}[b]{0.333\textwidth}
\centering

\includegraphics[width=\textwidth]{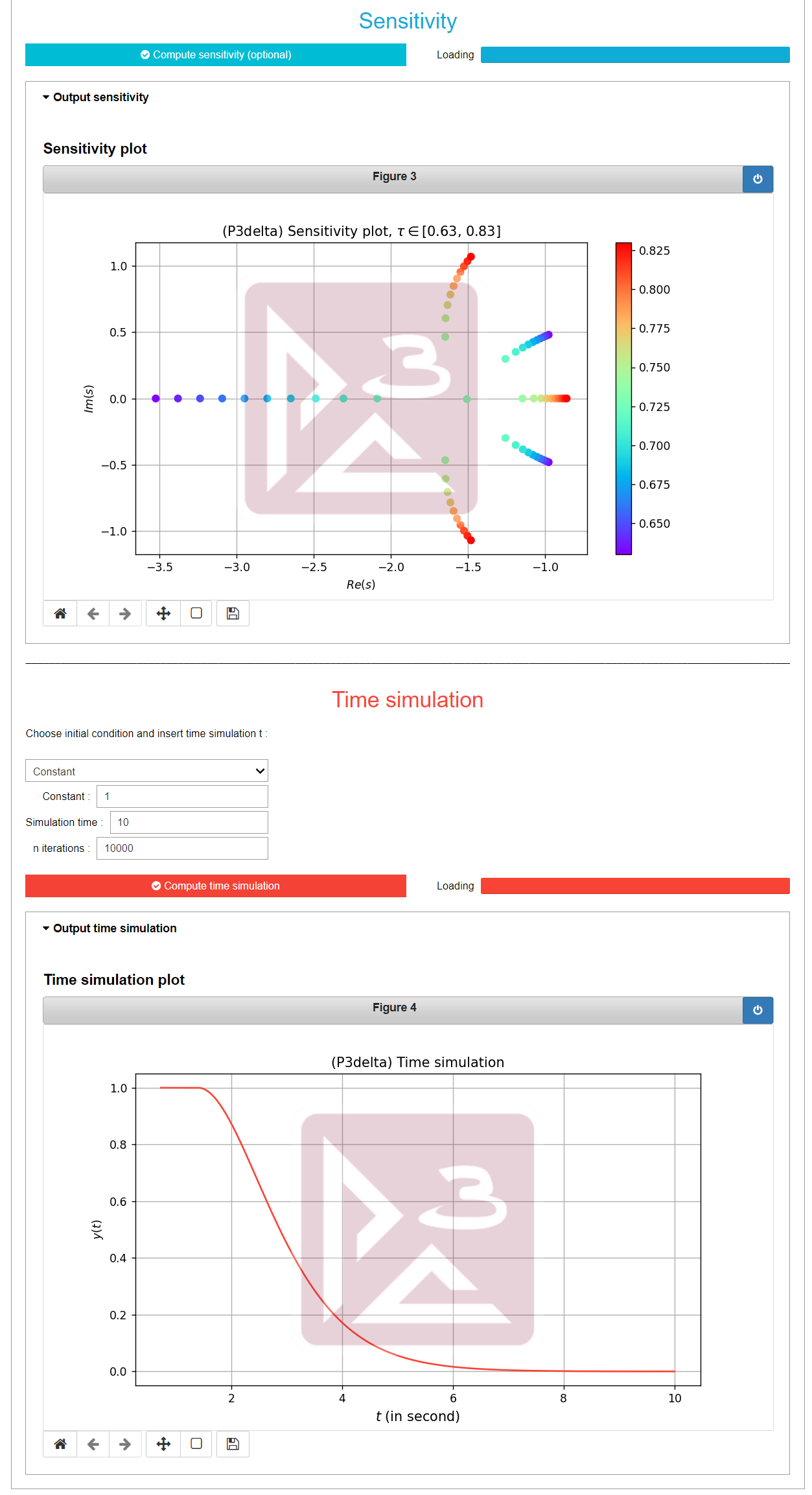}

(c)
\end{minipage}
\caption{{Example of use of P3$\delta$ Online in the ``Control-oriented MID'' mode. The user can easily select values of $n$ and $m$ using sliders and enter precise values for $s_0$ and $\tau$. The perturbed spectrum in (c) is computed with respect to variations in the delay value $\tau$. Computations are performed in real time and the user can switch between different figures and results. P3$\delta$ Online also generates a PDF report with data and figures.}}
\label{FigExplOnline}
\end{figure*}

In the ``Generic MID'' mode, the online version of P3$\delta$ returns the spectrum distribution as well as a normalized quasipolynomial which admits a root of multiplicity $n+m+1$ at the origin. In the ``Control-oriented MID'' mode, the online version of P3$\delta$ returns the admissibility region, a normalized quasipolynomial which admits a root of multiplicity $m+2$ at the origin, and an illustration of the bifurcation of the root of multiplicity $m+2$ with respect to variations of the value of the delay $\tau$. Illustrations of the online version of P3$\delta$ are provided in Figure~\ref{FigExplOnline}.

\section{Conclusion and planned developments}
The current version of P3$\delta$ exploits both MID and CRRID properties in its ``Generic'' mode and only the MID property in its the ``Control-oriented'' mode. Its main novelties with respect to its previous versions are the design for both retarded as well as neutral equations, the plot of the admissible assignment region in ``Control-oriented MID'' mode, and an online version of the software. In future developments, inspired from ``Control-oriented MID'' mode, other configurations for the spectrum distribution guaranteeing the dominancy of an assigned spectral value will be proposed, such as a ``Control-oriented CRRID'' mode.

\section*{Acknowledgments}
This work is partially supported by a public grant overseen by the French National Research Agency (ANR) as part of the ``Investissement d'Avenir'' program, through the iCODE project funded by the IDEX Paris-Saclay, ANR-11-IDEX0003-02. The authors also acknowledge the support of Institut Polytechnique des Sciences Avanc\'ees (IPSA).

The authors wish to acknowledge the work of the full P3$\delta$ development team, which, in addition to the authors, also include Mickael Alcaniz, Yoann Audet, Thomas Charbonnet, Honor\'e Curlier, Ayrton Hammoumou, Julien Huynh,  Pierre-Henry Poret, Achrafy Said Mohamed and Franck Sim. The development of P3$\delta$ was also made possible thanks to the work of the Cyb'Air Association.

\bibliography{p3delta}

\begin{thebibliography}{34}
\providecommand{\natexlab}[1]{#1}
\providecommand{\url}[1]{\texttt{#1}}
\providecommand{\urlprefix}{URL }
\expandafter\ifx\csname urlstyle\endcsname\relax
  \providecommand{\doi}[1]{doi:\discretionary{}{}{}#1}\else
  \providecommand{\doi}{doi:\discretionary{}{}{}\begingroup
  \urlstyle{rm}\Url}\fi

\bibitem[{Amrane et~al.(2018)Amrane, Bedouhene, Boussaada, and
  Niculescu}]{Amrane2018Qualitative}
Amrane, S., Bedouhene, F., Boussaada, I., and Niculescu, S.I. (2018).
\newblock On qualitative properties of low-degree quasipolynomials: further
  remarks on the spectral abscissa and rightmost-roots assignment.
\newblock \emph{Bull. Math. Soc. Sci. Math. Roumanie (N.S.)}, 61(109)(4),
  361--381.

\bibitem[{Atay(1999)}]{Atay1999Balancing}
Atay, F.M. (1999).
\newblock Balancing the inverted pendulum using position feedback.
\newblock \emph{Appl. Math. Lett.}, 12(5), 51--56.

\bibitem[{Avanessoff et~al.(2014)Avanessoff, Fioravanti, Bonnet, and
  Nguyen}]{Avanessoff2014Hinfty}
Avanessoff, D., Fioravanti, A.R., Bonnet, C., and Nguyen, L.H.V. (2014).
\newblock {$H_\infty$}-stability analysis of (fractional) delay systems of
  retarded and neutral type with the {M}atlab toolbox {YALTA}.
\newblock In \emph{Delay systems}, volume~1 of \emph{Adv. Delays Dyn.},
  285--297. Springer, Cham.

\bibitem[{Balogh et~al.(2021)Balogh, Boussaada, Insperger, and
  Niculescu}]{Balogh21}
Balogh, T., Boussaada, I., Insperger, T., and Niculescu, S.I. (2021).
\newblock Conditions for stabilizability of time-delay systems with real-rooted
  plant.
\newblock Submitted.

\bibitem[{Balogh et~al.(2020)Balogh, Insperger, Boussaada, and
  Niculescu}]{Balogh20}
Balogh, T., Insperger, T., Boussaada, I., and Niculescu, S.I. (2020).
\newblock Towards an {MID}-based delayed design for arbitrary-order dynamical
  systems with a mechanical application.
\newblock \emph{IFAC-PapersOnLine}, 53(2), 4375--4380.
\newblock 21th IFAC World Congress.

\bibitem[{Bedouhene et~al.(2020)Bedouhene, Boussaada, and
  Niculescu}]{BedouheneReal}
Bedouhene, F., Boussaada, I., and Niculescu, S.I. (2020).
\newblock Real spectral values coexistence and their effect on the stability of
  time-delay systems: Vandermonde matrices and exponential decay.
\newblock \emph{Comptes Rendus. Math\'ematique}, 358(9-10), 1011--1032.

\bibitem[{Benarab et~al.(2020)Benarab, Boussaada, Trabelsi, Mazanti, and
  Bonnet}]{Benarab2020MID}
Benarab, A., Boussaada, I., Trabelsi, K., Mazanti, G., and Bonnet, C. (2020).
\newblock The {MID} property for a second-order neutral time-delay differential
  equation.
\newblock In \emph{2020 24th International Conference on System Theory, Control
  and Computing (ICSTCC)}, 202--207.

\bibitem[{Berenstein and Gay(1995)}]{Berenstein1995Complex}
Berenstein, C.A. and Gay, R. (1995).
\newblock \emph{Complex analysis and special topics in harmonic analysis}.
\newblock Springer-Verlag, New York.

\bibitem[{Boussaada et~al.(2021)Boussaada, Mazanti, and
  Niculescu}]{BoussaadaGeneric}
Boussaada, I., Mazanti, G., and Niculescu, S.I. (2021).
\newblock The generic multiplicity-induced-dominancy property from retarded to
  neutral delay-differential equations: When delay-systems characteristics meet
  the zeros of {K}ummer functions.
\newblock Submitted.

\bibitem[{Boussaada et~al.(2020{\natexlab{a}})Boussaada, Mazanti, Niculescu,
  Huynh, Sim, and Thomas}]{BoussaadaP3D}
Boussaada, I., Mazanti, G., Niculescu, S.I., Huynh, J., Sim, F., and Thomas, M.
  (2020{\natexlab{a}}).
\newblock Partial pole placement via delay action: A {P}ython software for
  delayed feedback stabilizing design.
\newblock In \emph{2020 24th International Conference on System Theory, Control
  and Computing (ICSTCC)}, 196--201.

\bibitem[{Boussaada and Niculescu(2018)}]{Boussaada2018Dominancy}
Boussaada, I. and Niculescu, S.I. (2018).
\newblock On the dominancy of multiple spectral values for time-delay systems
  with applications.
\newblock \emph{IFAC-PapersOnLine}, 51(14), 55--60.
\newblock 14th IFAC Workshop on Time Delay Systems.

\bibitem[{Boussaada et~al.(2020{\natexlab{b}})Boussaada, Niculescu, El-Ati,
  P\'{e}rez-Ramos, and Trabelsi}]{Boussaada2020Multiplicity}
Boussaada, I., Niculescu, S.I., El-Ati, A., P\'{e}rez-Ramos, R., and Trabelsi,
  K. (2020{\natexlab{b}}).
\newblock Multiplicity-induced-dominancy in parametric second-order delay
  differential equations: {A}nalysis and application in control design.
\newblock \emph{ESAIM Control Optim. Calc. Var.}, 26, Paper No. 57.

\bibitem[{Boussaada et~al.(2018)Boussaada, Tliba, Niculescu, \"{U}nal, and
  Vyhl\'{\i}dal}]{Boussaada2018Further}
Boussaada, I., Tliba, S., Niculescu, S.I., \"{U}nal, H.U., and Vyhl\'{\i}dal,
  T. (2018).
\newblock Further remarks on the effect of multiple spectral values on the
  dynamics of time-delay systems. {A}pplication to the control of a mechanical
  system.
\newblock \emph{Linear Algebra Appl.}, 542, 589--604.

\bibitem[{Breda et~al.(2009)Breda, Maset, and Vermiglio}]{Breda2009Trace}
Breda, D., Maset, S., and Vermiglio, R. (2009).
\newblock T{RACE}-{DDE}: a tool for robust analysis and characteristic
  equations for delay differential equations.
\newblock In \emph{Topics in time delay systems}, volume 388 of \emph{Lect.
  Notes Control Inf. Sci.}, 145--155. Springer, Berlin.

\bibitem[{Engelborghs et~al.(2002)Engelborghs, Luzyanina, and
  Roose}]{Engelborghs2002Numerical}
Engelborghs, K., Luzyanina, T., and Roose, D. (2002).
\newblock Numerical bifurcation analysis of delay differential equations using
  {DDE}-{BIFTOOL}.
\newblock \emph{ACM Trans. Math. Software}, 28(1), 1--21.

\bibitem[{Hale and Verduyn~Lunel(1993)}]{Hale1993Introduction}
Hale, J.K. and Verduyn~Lunel, S.M. (1993).
\newblock \emph{Introduction to functional differential equations}.
\newblock Springer-Verlag, New York.

\bibitem[{Irofti et~al.(2016)Irofti, Boussaada, and
  Niculescu}]{Irofti2016Codimension}
Irofti, D.A., Boussaada, I., and Niculescu, S.I. (2016).
\newblock On the codimension of the singularity at the origin for networked
  delay systems.
\newblock In A.~Seuret, L.~Hetel, J.~Daafouz, and K.H. Johansson (eds.),
  \emph{Delays and Networked Control Systems}, 3--15. Springer International
  Publishing.

\bibitem[{Ma et~al.(2020)Ma, Boussaada, Bonnet, Niculescu, and Chen}]{Ma}
Ma, D., Boussaada, I., Bonnet, C., Niculescu, S.I., and Chen, J. (2020).
\newblock {Multiplicity-Induced-Dominancy extended to neutral delay equations:
  Towards a systematic PID tuning based on Rightmost root assignment}.
\newblock In \emph{{ACC 2020 - American Control Conference}}. Denver, United
  States.

\bibitem[{Manitius and Olbrot(1979)}]{Manitius1979Finite}
Manitius, A.Z. and Olbrot, A.W. (1979).
\newblock Finite spectrum assignment problem for systems with delays.
\newblock \emph{IEEE Trans. Automat. Control}, 24(4), 541--553.

\bibitem[{Mazanti et~al.(2020{\natexlab{a}})Mazanti, Boussaada, and
  Niculescu}]{Mazanti2020Qualitative}
Mazanti, G., Boussaada, I., and Niculescu, S.I. (2020{\natexlab{a}}).
\newblock On qualitative properties of single-delay linear retarded
  differential equations: Characteristic roots of maximal multiplicity are
  necessarily dominant.
\newblock \emph{IFAC-PapersOnLine}, 53(2), 4345--4350.
\newblock 21st IFAC World Congress.

\bibitem[{Mazanti et~al.(2021{\natexlab{a}})Mazanti, Boussaada, and
  Niculescu}]{MazantiMultiplicity}
Mazanti, G., Boussaada, I., and Niculescu, S.I. (2021{\natexlab{a}}).
\newblock Multiplicity-induced-dominancy for delay-differential equations of
  retarded type.
\newblock \emph{J. Differential Equations}, 286, 84--118.

\bibitem[{Mazanti et~al.(2021{\natexlab{b}})Mazanti, Boussaada, Niculescu, and
  Chitour}]{MazantiEffects}
Mazanti, G., Boussaada, I., Niculescu, S.I., and Chitour, Y.
  (2021{\natexlab{b}}).
\newblock Effects of roots of maximal multiplicity on the stability of some
  classes of delay differential-algebraic systems: the lossless propagation
  case.
\newblock \emph{IFAC-PapersOnLine}, 54(9), 764--769.
\newblock 24th International Symposium on Mathematical Theory of Networks and
  Systems (MTNS 2020).

\bibitem[{Mazanti et~al.(2020{\natexlab{b}})Mazanti, Boussaada, Niculescu, and
  Vyhl{\'{\i}}dal}]{Mazanti2020Spectral}
Mazanti, G., Boussaada, I., Niculescu, S.I., and Vyhl{\'{\i}}dal, T.
  (2020{\natexlab{b}}).
\newblock Spectral dominance of complex roots for single-delay linear
  equations.
\newblock \emph{IFAC-PapersOnLine}, 53(2), 4357--4362.
\newblock 21st IFAC World Congress.

\bibitem[{Michiels et~al.(2002)Michiels, Engelborghs, Vansevenant, and
  Roose}]{Michiels2002}
Michiels, W., Engelborghs, K., Vansevenant, P., and Roose, D. (2002).
\newblock Continuous pole placement for delay equations.
\newblock \emph{Automatica}, 38(5), 747--761.

\bibitem[{Michiels and Niculescu(2014)}]{Michiels2014Stability}
Michiels, W. and Niculescu, S.I. (2014).
\newblock \emph{Stability, control, and computation for time-delay systems: An
  eigenvalue-based approach}.
\newblock SIAM, Philadelphia, PA, second edition.

\bibitem[{Niculescu et~al.(2010)Niculescu, Michiels, Gu, and
  Abdallah}]{Niculescu2010Delay}
Niculescu, S.I., Michiels, W., Gu, K., and Abdallah, C.T. (2010).
\newblock Delay effects on output feedback control of dynamical systems.
\newblock In F.M. Atay (ed.), \emph{Complex time-delay systems}, 63--84.
  Springer, Berlin.

\bibitem[{Parini(2018--)}]{cxroots}
Parini, R. (2018--).
\newblock {cxroots: A Python module to find all the roots of a complex analytic
  function within a given contour}.
\newblock \urlprefix\url{https://github.com/rparini/cxroots}.

\bibitem[{{P}roject {J}upyter et~al.(2018){P}roject {J}upyter, {M}atthias
  {B}ussonnier, {J}essica {F}orde, {J}eremy {F}reeman, {B}rian {G}ranger, {T}im
  {H}ead, {C}hris {H}oldgraf, {K}yle {K}elley, {G}ladys {N}alvarte, {A}ndrew
  {O}sheroff, {P}acer, {Y}uvi {P}anda, {F}ernando {P}erez, {B}enjamin~{R}agan
  {K}elley, and {C}arol {W}illing}]{project_jupyter-proc-scipy-2018}
{P}roject {J}upyter, {M}atthias {B}ussonnier, {J}essica {F}orde, {J}eremy
  {F}reeman, {B}rian {G}ranger, {T}im {H}ead, {C}hris {H}oldgraf, {K}yle
  {K}elley, {G}ladys {N}alvarte, {A}ndrew {O}sheroff, {P}acer, M., {Y}uvi
  {P}anda, {F}ernando {P}erez, {B}enjamin~{R}agan {K}elley, and {C}arol
  {W}illing (2018).
\newblock {B}inder 2.0 - {R}eproducible, interactive, sharable environments for
  science at scale.
\newblock In {F}atih {A}kici, {D}avid {L}ippa, {D}illon {N}iederhut, and
  M.~{P}acer (eds.), \emph{{P}roceedings of the 17th {P}ython in {S}cience
  {C}onference}, 113--120.

\bibitem[{Ram et~al.(2011)Ram, Mottershead, and Tehrani}]{RAMLAA}
Ram, Y., Mottershead, J., and Tehrani, M. (2011).
\newblock Partial pole placement with time delay in structures using the
  receptance and the system matrices.
\newblock \emph{Linear Algebra and its Applications}, 434(7), 1689--1696.

\bibitem[{St\'{e}p\'{a}n(1989)}]{Stepan1989Retarded}
St\'{e}p\'{a}n, G. (1989).
\newblock \emph{Retarded dynamical systems: stability and characteristic
  functions}, volume 210 of \emph{Pitman Research Notes in Mathematics Series}.
\newblock Longman Scientific \& Technical, Harlow.

\bibitem[{Suh and Bien(1979)}]{Suh1979Proportional}
Suh, I.H. and Bien, Z. (1979).
\newblock Proportional minus delay controller.
\newblock \emph{IEEE Trans. Automat. Control}, 24(2), 370--372.

\bibitem[{Tallman and Smith(1958)}]{Tallman1958Analog}
Tallman, G.H. and Smith, O.J.M. (1958).
\newblock Analog study of dead-beat posicast control.
\newblock \emph{IRE Transactions on Automatic Control}, 4(1), 14--21.

\bibitem[{Vyhl\'{\i}dal and Z\'{\i}tek(2014)}]{Vyhlidal2014QPmR}
Vyhl\'{\i}dal, T. and Z\'{\i}tek, P. (2014).
\newblock Q{P}m{R}---quasi-polynomial root-finder: algorithm update and
  examples.
\newblock In \emph{Delay systems}, volume~1 of \emph{Adv. Delays Dyn.},
  299--312. Springer, Cham.

\bibitem[{Wielonsky(2001)}]{Wielonsky2001Rolle}
Wielonsky, F. (2001).
\newblock A {R}olle's theorem for real exponential polynomials in the complex
  domain.
\newblock \emph{J. Math. Pures Appl. (9)}, 80(4), 389--408.

\end{thebibliography}

\end{document}